\documentclass{article}

\usepackage{amssymb,amsmath,amsthm}

\newtheorem{thm}{Theorem}[section]
\theoremstyle{definition}
\newtheorem{exmp}[thm]{Example}
\newtheorem{defn}[thm]{Definition}

\theoremstyle{remark}
\newtheorem{rem}[thm]{Remark}

\begin{document}

\title{A Formulation of Noether's Theorem
for Fractional Problems
of the Calculus of Variations\footnote{Research report
CM06/I-04, University of Aveiro. Accepted for publication
in the \emph{Journal of Mathematical Analysis and Applications.}}}

\author{Gast\~{a}o S. F. Frederico\\
\texttt{gfrederico@mat.ua.pt} \and
Delfim F. M. Torres\\
\texttt{delfim@mat.ua.pt}}

\date{Department of Mathematics\\
University of Aveiro\\
3810-193 Aveiro, Portugal}

\maketitle


\begin{abstract}
Fractional (or non-integer) differentiation is an important concept
both from theoretical and applicational points of view. The study
of problems of the calculus of variations with fractional
derivatives is a rather recent subject, the main result being the
fractional necessary optimality condition of Euler-Lagrange
obtained in 2002. Here we use the notion of Euler-Lagrange
fractional extremal to prove a Noether-type theorem. For that we
propose a generalization of the classical concept of conservation
law, introducing an appropriate fractional operator.

\medskip

\noindent \textbf{Keywords:} calculus of variations;
fractional derivatives; Noether's theorem.

\medskip

\noindent \textbf{Mathematics Subject Classification 2000:} 49K05; 26A33.

\end{abstract}


\section{Introduction}

The notion of conservation law -- first integral
of the Euler-Lagrange equations -- is well-known in
Physics. One of the most important conservation laws
is the integral of energy, discovered by Leonhard Euler in 1744:
when a Lagrangian $L(q,\dot{q})$ corresponds to a system of
conservative points, then
\begin{equation}
\label{eq:R}
\begin{gathered}
-L(q,\dot{q})+\partial_{2} L(q,\dot{q})\cdot\dot{q}
\equiv \text{constant}
\end{gathered}
\end{equation}
holds along all the solutions of the Euler-Lagrange equations
(along the extremals of the autonomous variational problem),
where $\partial_{2} L(\cdot,\cdot)$ denote the partial derivative
of function $L(\cdot,\cdot)$ with respect to its second argument.
Many other examples appear in modern physics: in classic, quantum, and
optical mechanics; in the theory of relativity; etc. For instance,
in classic mechanics, beside the conservation of energy
\eqref{eq:R}, it may occur conservation of momentum or angular
momentum. These conservation laws are very important: they can be
used to reduce the order of the Euler-Lagrange differential
equations, thus simplifying the resolution of the problems.

In 1918 Emmy Noether proved a general theorem of the calculus
of variations, that permits to obtain, from the existence
of variational symmetries, all the conservation laws
that appear in applications. In the last decades,
Noether's principle has been formulated
in various contexts (see \cite{delfimPortMath04,delfimCPAA04}
and references therein). In this work we generalize Noether's
theorem for problems having fractional derivatives.

Fractional differentiation plays nowadays an important role in
various seemingly diverse and widespread fields of science and
engineering: physics (classic and quantum mechanics,
thermodynamics, optics, etc), chemistry, biology, economics,
geology, astrophysics, probability and statistics,
signal and image processing, dynamics of earthquakes,
control theory, and so on \cite{Agrawal:2004b,CD:Hilfer:2000,Kilbas,CD:Klimek:2002}.
Its origin goes back more than 300 years, when in 1695 L'Hopital asked Leibniz the
meaning of $\frac{d^{n}y}{dx^{n}}$ for $n=\frac{1}{2}$. After
that, many famous mathematicians, like J.~Fourier,
N.~H.~Abel, J.~Liouville, B.~Riemann, among others, contributed to
the development of the Fractional Calculus
\cite{CD:Hilfer:2000,CD:MilRos:1993,CD:SaKiMa:1993}.

F.~Riewe \cite{CD:Riewe:1996,CD:Riewe:1997} obtained a version of
the Euler-Lagrange equations for problems of the Calculus of
Variations with fractional derivatives, that combines the
conservative and non-conservative cases. More recently, O.~Agrawal
proved a formulation for variational problems with right and left
fractional derivatives in the Riemann-Liouville sense
\cite{CD:Agrawal:2002}. Then these Euler-Lagrange equations were
used by D.~Baleanu and T.~Avkar to investigate problems with
Lagrangians which are linear on the velocities
\cite{CD:BalAv:2004}. Here we use the results of
\cite{CD:Agrawal:2002} to generalize Noether's theorem for the
more general context of the Fractional Calculus of Variations.

The paper is organized in the following way. In
Section~\ref{sec:fdRL} we recall the notions of right and left
Riemann-Liouville fractional derivatives, that are needed for
formulating the fractional problem of the calculus of
variations. There are many different ways to approach classical
Noether's theorem. In Section~\ref{sec:cNT} we review the only
proof that we are able to extend, with success, to the fractional
context. The method is based on a two-step procedure: it starts
with an invariance notion of the integral functional under a
one-parameter infinitesimal group of transformations, without
changing the time variable; then it proceeds with a
time-reparameterization to obtain Noether's theorem in
general form. The intended fractional Noether's theorem is
formulated and proved in Section~\ref{sec:MR}. Two illustrative
examples of application of our main result are given in Section~\ref{sec:Ex}.
We finish with Section~\ref{sec:Conc} of conclusions and some open questions.


\section{Riemann-Liouville fractional derivatives}
\label{sec:fdRL}

In this section we collect the definitions of right and left
Riemann-Liouville fractional derivatives and their main properties
\cite{CD:Agrawal:2002,CD:MilRos:1993,CD:SaKiMa:1993}.

\begin{defn}[Riemann-Liouville fractional derivatives]
Let $f$ be a continuous and integrable function in the interval
$[a,b]$. For all $t \in [a,b]$, the left Riemann-Liouville
fractional derivative $_aD_t^\alpha f(t)$, and the right
Riemann-Liouville fractional derivative $_tD_b^\alpha f(t)$, of
order $\alpha$, are defined in the following way:
\begin{gather}
_aD_t^\alpha f(t)
= \frac{1}{\Gamma(n-\alpha)}\left(\frac{d}{dt}\right)^{n}
\int_a^t (t-\theta)^{n-\alpha-1}f(\theta)d\theta \, , \label{eq:DFRLE} \\
_tD_b^\alpha f(t)
=\frac{1}{\Gamma(n-\alpha)}\left(-\frac{d}{dt}\right)^{n}
\int_t^b(\theta - t)^{n-\alpha-1}f(\theta)d\theta \, , \label{eq:DFRLD}
\end{gather}
where $n \in \mathbb{N}$, $n-1 \leq \alpha < n$, and $\Gamma$ is
the Euler gamma function.
\end{defn}

\begin{rem}
If $\alpha$ is an integer, then from \eqref{eq:DFRLE}
and \eqref{eq:DFRLD} one obtains the standard derivatives, that is,
\begin{equation*}
\label{eq:DU}
_aD_t^\alpha f(t) = \left(\frac{d}{dt}\right)^\alpha f(t) \, , \quad
_tD_b^\alpha f(t) = \left(-\frac{d}{dt}\right)^\alpha f(t) \, .
\end{equation*}
\end{rem}

\begin{thm} Let $f$ and $g$ be two continuous functions
on $[a,b]$. Then, for all $t \in [a,b]$, the following properties hold:
\begin{enumerate}
\item for $p>0$, $_aD_t^p\left(f(t)+g(t)\right)
= {_aD_t^p}f(t)+{_aD_t^p}g(t)$;
\item for $p \geq q \geq 0$,
$_aD_t^p\left(_aD_t^{-q} f(t)\right) = {_aD_t^{p-q}}f(t)$;
\item for $p>0$, $_aD_t^p\left(_aD_t^{-p} f(t)\right) = f(t)$
(fundamental property of the Riemann-Liouville fractional derivatives);
\item for $p>0$, $\int_a^b \left(_aD_t^p f(t)\right)g(t)dt
= \int_a^b f(t) \, _tD_b^p g(t)dt$ (we are assuming that $_aD_t^p f(t)$
and $_tD_b^p g(t)$ exist at every point $t \in [a,b]$ and are continuous).
\end{enumerate}
\end{thm}

\begin{rem}
In general, the fractional derivative of a constant is not
equal to zero.
\end{rem}

\begin{rem}
The fractional derivative of order $p>0$
of function $(t-a)^\upsilon$, $\upsilon>-1$,
is given by
\begin{equation*}
_aD_t^p(t-a)^\upsilon
= \frac{\Gamma(\upsilon+1)}{\Gamma(-p+\upsilon+1)}(t-a)^{\upsilon-p} \, .
\end{equation*}
\end{rem}

\begin{rem}
In the literature, when one reads ``Riemann-Liouville fractional derivative'',
one usually means the ``left Riemann-Liouville fractional derivative''.
In Physics, if $t$ denotes the time-variable, the
right Riemann-Liouville fractional derivative of $f(t)$
is interpreted as a future state of the process $f(t)$. For this
reason, the right-derivative is usually
neglected in applications, when the present state of the process does
not depend on the results of the future development. From a mathematical
point of view, both derivatives appear naturally in the fractional calculus of
variations \cite{CD:Agrawal:2002}.
\end{rem}

We refer the reader interested in additional background
on fractional theory, to the comprehensive book
of Samko et al. \cite{CD:SaKiMa:1993}.


\section{Review of the classical Noether's theorem}
\label{sec:cNT}

There exist several ways to prove the classical theorem of Emmy
Noether. In this section we review one of those proofs
\cite{CD:Jost:1998}. The proof is done in two steps: we begin by
proving Noether's theorem without transformation of the time
(without transformation of the independent variable); then, using
a technique of time-reparameterization, we obtain Noether's
theorem in its general form. This technique is not so popular
while proving Noether's theorem, but it turns out to be the only
one we succeed to extend for the more general context of the
Fractional Calculus of Variations.

We begin by formulating the fundamental problem
of the calculus of variations:
\begin{equation}
\label{P}
I[q(\cdot)] = \int_a^b L\left(t,q(t),\dot{q}(t)\right) dt
\longrightarrow \min  \tag{P}
\end{equation}
under the boundary conditions $q(a)=q_{a}$ and $q(b)=q_{b}$,
and where $\dot{q} = \frac{dq}{dt}$. The Lagrangian
$L :[a,b] \times \mathbb{R}^{n} \times \mathbb{R}^{n} \rightarrow \mathbb{R}$
is assumed to be a $C^{2}$-function with respect to all its arguments.

\begin{defn}[invariance without transforming the time]
\label{def:inv}
Functional \eqref{P} is said to be invariant under a $\varepsilon$-parameter
group of infinitesimal transformations
\begin{equation}
\label{eq:tinf}
\bar{q}(t)= q(t) + \varepsilon\xi(t,q) + o(\varepsilon)
\end{equation}
if, and only if,
\begin{equation}
\label{eq:inv}
\int_{t_{a}}^{t_{b}} L\left(t,q(t),\dot{q}(t)\right)dt
= \int_{t_{a}}^{t_{b}}L\left(t,\bar{q}(t),\dot{\bar{q}}(t)\right)dt
\end{equation}
for any subinterval $[{t_{a}},{t_{b}}] \subseteq [a,b]$.
\end{defn}

Along the work, we denote by $\partial_{i}L$ the partial derivative
of $L$ with respect to its $i$-th argument.

\begin{thm}[necessary condition of invariance]
\label{theo:cnsi}
If functional \eqref{P} is invariant under transformations
\eqref{eq:tinf}, then
\begin{equation}
\label{eq:cnsi}
\partial_{2} L\left(t,q,\dot{q}\right)\cdot\xi
+\partial_{3} L\left(t,q,\dot{q}\right)\cdot\dot{\xi} = 0 \, .
\end{equation}
\end{thm}

\begin{proof}
Equation \eqref{eq:inv} is equivalent to
\begin{equation}
\label{eq:inv1}
L\left(t,q,\dot{q}\right) =
L(t,q+\varepsilon\xi+o(\varepsilon),\dot{q}
+ \varepsilon\dot{\xi} + o(\varepsilon)) \, .
\end{equation}
Differentiating both sides of equation
\eqref{eq:inv1} with respect to $\varepsilon$,
then substituting $\varepsilon=0$, we obtain equality \eqref{eq:cnsi}.
\end{proof}

\begin{defn}[conserved quantity]
\label{def:leicon}
Quantity $C(t,q(t),\dot{q}(t))$ is said to be \emph{conserved}
if, and only if, $\frac{d}{dt}C(t,q(t),\dot{q}(t))=0$ along all the
solutions of the Euler-Lagrange equations
\begin{equation}
\label{eq:el}
\frac{d}{{dt}}\partial_{3} L\left(t,q,\dot{q}\right)
= \partial_{2} L\left(t,q,\dot{q}\right) \, .
\end{equation}
\end{defn}

\begin{thm}[Noether's theorem without transforming time]
\label{theo:tnaut}
If functional \eqref{P} is invariant under the one-parameter
group of transformations \eqref{eq:tinf}, then
\begin{equation}
\label{eq:tnau}
C(t,q,\dot{q}) = \partial_{3} L\left(t,q,\dot{q}\right)\cdot\xi(t,q)
\end{equation}
is conserved.
\end{thm}

\begin{proof}
Using the Euler-Lagrange equations \eqref{eq:el}
and the necessary condition of invariance \eqref{eq:cnsi}, we obtain:
\begin{equation*}
\begin{split}
\frac{d}{dt}&\left(\partial_{3} L\left(t,q,\dot{q}\right) \cdot \xi(t,q)\right) \\
&= \frac{d}{{dt}}\partial_{3} L\left(t,q,\dot{q}\right) \cdot \xi(t,q)
+ \partial_{3} L\left(t,q,\dot{q}\right) \cdot \dot{\xi}(t,q) \\
&=  \partial_{2} L\left(t,q,\dot{q}\right) \cdot \xi(t,q)
+ \partial_{3} L\left(t,q,\dot{q}\right) \cdot \dot{\xi}(t,q) \\
&= 0 \, .
\end{split}
\end{equation*}
\end{proof}

\begin{rem}
In classical mechanics, $\partial_{3} L\left(t,q,\dot{q}\right)$
is interpreted as the generalized momentum.
\end{rem}

\begin{defn}[invariance of (\ref{P})]
\label{def:inva}
Functional $(P)$ is said to be invariant under the
one-parameter group of infinitesimal transformations
\begin{equation}
\label{eq:tinf2}
\begin{cases}
\bar{t} = t + \varepsilon\tau(t,q) + o(\varepsilon) \, ,\\
\bar{q}(t) = q(t) + \varepsilon\xi(t,q) + o(\varepsilon) \, ,\\
\end{cases}
\end{equation}
if, and only if,
\begin{equation*}
\label{eq:inv2}
\int_{t_{a}}^{t_{b}} L\left(t,q(t),\dot{q}(t)\right)dt
= \int_{\bar{t}(t_a)}^{\bar{t}(t_b)}
L\left(\bar{t},\bar{q}(\bar{t}),\dot{\bar{q}}(\bar{t})\right)d\bar{t}
\end{equation*}
for any subinterval $[{t_{a}},{t_{b}}] \subseteq [a,b]$.
\end{defn}

\begin{thm}[Noether's theorem]
\label{theo:tnoe}
If functional \eqref{P} is invariant, in the sense of
Definition~\ref{def:inva}, then
\begin{multline}
\label{eq:TeoNet}
C(t,q,\dot{q}) =
\partial_{3} L\left(t,q,\dot{q}\right)\cdot\xi(t,q)
+ \left( L(t,q,\dot{q}) - \partial_{3} L\left(t,q,\dot{q}\right)
\cdot \dot{q} \right) \tau(t,q)
\end{multline}
is conserved.
\end{thm}

\begin{proof}
Every non-autonomous problem \eqref{P}
is equivalent to an autonomous one,
considering $t$ as a dependent variable. For that
we consider a Lipschitzian one-to-one transformation
\begin{equation*}
[a,b]\ni t \longmapsto \sigma \in [\sigma_{a},\sigma_{b}]
\end{equation*}
such that
\begin{equation*}
\begin{split}
I\left[q(\cdot)\right] &= \int_a^b L\left(t,q(t),\dot{q}(t)\right)dt \\
&= \int_{\sigma_{a}}^{\sigma_{b}} L\left(t(\sigma),q(t(\sigma)),
\frac{\frac{dq(t(\sigma))}{d\sigma}}{\frac{dt(\sigma)}{d\sigma}}\right)
\frac{dt(\sigma)}{d\sigma}d\sigma \\
&= \int_{\sigma_{a}}^{\sigma_{b}}
L\left(t(\sigma),q(t(\sigma)),\frac{q_{\sigma}^{'}}{t_{\sigma}^{'}}\right)t_{\sigma}^{'}d\sigma\\
&\doteq \int_{\sigma_{a}}^{\sigma_{b}}
\bar{L}\left(t(\sigma),q(t(\sigma)),t_{\sigma}^{'},q_{\sigma}^{'}\right)d\sigma \\
&\doteq \bar{I}\left[t(\cdot),q(t(\cdot))\right] \, ,
\end{split}
\end{equation*}
where $t(\sigma_{a})=a$, $t(\sigma_{b})=b$,
$t_{\sigma}^{'}=\frac{dt(\sigma)}{d\sigma}$, and
$q_{\sigma}^{'}=\frac{dq(t(\sigma))}{d\sigma}$.
If functional $I[q(\cdot)]$ is invariant in the sense
of Definition~\ref{def:inva}, then functional
$\bar{I}[t(\cdot),q(t(\cdot))]$ is invariant in the sense
of Definition~\ref{def:inv}. Applying
Theorem~\ref{theo:tnaut}, we obtain that
\begin{equation}
\label{eq:dtn}
C\left(t,q,t_{\sigma}^{'},q_{\sigma}^{'}\right)
= \partial_{4} \bar{L}\cdot\xi + \partial_{3} \bar{L}\tau
\end{equation}
is a conserved quantity. Since
\begin{equation}
\label{eq:dtn1}
\begin{split}
\partial_{4} \bar{L} &= \partial_{3} L\left(t,q,\dot{q}\right) \, ,\\
\partial_{3} \bar{L} &= -\partial_{3} L\left(t,q,\dot{q}\right)
\cdot \frac{q_{\sigma}^{'}}{t_{\sigma}^{'}}+L(t,q,\dot{q}) \\
&= L(t,q,\dot{q})-\partial_{3} L\left(t,q,\dot{q}\right)\cdot\dot{q} \, ,
\end{split}
\end{equation}
substituting \eqref{eq:dtn1} into \eqref{eq:dtn},
we arrive to the intended conclusion \eqref{eq:TeoNet}.
\end{proof}


\section{Main Results}
\label{sec:MR}

In 2002 \cite{CD:Agrawal:2002}, a formulation of the
Euler-Lagrange equations was given for problems of the calculus of
variations with fractional derivatives. In this section we prove a
Noether's theorem for the fractional Euler-Lagrange extremals.

The fundamental functional of the fractional calculus of variations
is defined as follows:
\begin{gather}
\label{Pf}
I[q(\cdot)] = \int_a^b L\left(t,q(t),{_aD_t^\alpha} q(t), {_tD_b^\beta} q(t)\right) dt
\longrightarrow \min \, , \tag{$P_f$}
\end{gather}
where the Lagrangian
$L :[a,b] \times \mathbb{R}^{n} \times
\mathbb{R}^{n} \times \mathbb{R}^{n} \rightarrow \mathbb{R}$
is a $C^{2}$ function with respect to all its arguments, and
$0 < \alpha,\beta \leq 1$.

\begin{rem}
In the case $\alpha = \beta = 1$, problem \eqref{Pf} is reduced to
problem \eqref{P}:
\begin{equation*}
I[q(\cdot)] = \int_a^b \mathcal{L}\left(t,q(t),\dot{q}(t)\right) dt
\longrightarrow \min
\end{equation*}
with
\begin{equation}
\label{eq:mathcalL}
\mathcal{L}\left(t,q,\dot{q}\right) = L(t,q,\dot{q},-\dot{q}) \, .
\end{equation}
\end{rem}

Theorem~\ref{Thm:FractELeq} summarizes
the main result of \cite{CD:Agrawal:2002}.

\begin{thm}[\cite{CD:Agrawal:2002}]
\label{Thm:FractELeq}
If $q$ is a minimizer of problem \eqref{Pf}, then
it satisfies the \emph{fractional Euler-Lagrange equations}:
\begin{multline}
\label{eq:eldf}
\partial_{2} L\left(t,q,{_aD_t^\alpha q},{_tD_b^\beta q}\right)
+ {_tD_b^\alpha}\partial_{3} L\left(t,q,{_aD_t^\alpha q},{_tD_b^\beta q}\right) \\
+ {_aD_t^\beta}\partial_{4} L\left(t,q,{_aD_t^\alpha q},{_tD_b^\beta q}\right) = 0 \, .
\end{multline}
\end{thm}

\begin{defn}[\textrm{cf.} Definition~\ref{def:inv}]
\label{def:inv1:MR}
We say that functional \eqref{Pf} is
invariant under the transformation \eqref{eq:tinf} if, and only if,
\begin{multline}
\label{eq:invdf}
\int_{t_{a}}^{t_{b}} L\left(t,q(t),{_aD_t^\alpha q(t)},{_tD_b^\beta q(t)}\right) dt
= \int_{t_{a}}^{t_{b}} L\left(t,\bar{q}(t),{_aD_t^\alpha
\bar{q}(t)},{_tD_b^\beta \bar{q}(t)}\right) dt
\end{multline}
for any subinterval $[{t_{a}},{t_{b}}] \subseteq [a,b]$.
\end{defn}

The next theorem establishes a necessary condition of
invariance, of extreme importance for our objectives.

\begin{thm}[\textrm{cf.} Theorem~\ref{theo:cnsi}]
If functional \eqref{Pf} is invariant under
transformations \eqref{eq:tinf}, then
\begin{multline}
\label{eq:cnsidf}
\partial_{2} L\left(t,q,{_aD_t^\alpha q},{_tD_b^\beta q}\right) \cdot \xi(t,q)
+ \partial_{3} L\left(t,q,{_aD_t^\alpha q},{_tD_b^\beta q}\right)
\cdot {_aD_t^\alpha \xi(t,q)} \\
+ \partial_{4} L\left(t,q,{_aD_t^\alpha q},{_tD_b^\beta q}\right)
\cdot {_tD_b^\beta} \xi(t,q) = 0 \, .
\end{multline}
\end{thm}

\begin{rem}
In the particular case $\alpha = \beta =1$, we obtain from
\eqref{eq:cnsidf} the necessary condition \eqref{eq:cnsi} applied to
$\mathcal{L}$ \eqref{eq:mathcalL}.
\end{rem}

\begin{proof}
Having in mind that condition \eqref{eq:invdf}
is valid for any subinterval $[{t_{a}},{t_{b}}] \subseteq [a,b]$,
we can get rid off the integral signs in \eqref{eq:invdf}.
Differentiating this condition
with respect to $\varepsilon$, substituting $\varepsilon=0$,
and using the definitions and properties
of the Riemann-Liouville fractional derivatives given in
Section~\ref{sec:fdRL}, we arrive to
\begin{multline}
\label{eq:SP}
0 = \partial_{2} L\left(t,q,{_aD_t^\alpha q},{_tD_b^\beta q}\right)\cdot\xi(t,q) \\
+ \partial_{3} L\left(t,q,{_aD_t^\alpha q},{_tD_b^\beta q}\right)
\cdot\frac{d}{d\varepsilon} \left[\frac{1}{\Gamma(n-\alpha)}
\left(\frac{d}{dt}\right)^{n}\int_a^t (t-\theta)^{n-\alpha-1}q(\theta)d\theta\right. \\
+\left.\frac{\varepsilon}{\Gamma(n-\alpha)}\left(\frac{d}{dt}\right)^{n}\int_a^t
(t-\theta)^{n-\alpha-1}\xi(\theta,q)d\theta\right]_{\varepsilon=0}  \\
+ \partial_{4} L\left(t,q,{_aD_t^\alpha q},{_tD_b^\beta q}\right)
\cdot\frac{d}{d\varepsilon}\left[\frac{1}{\Gamma(n-\beta)}\left(-\frac{d}{dt}\right)^{n}
\int_t^b (\theta - t)^{n-\beta-1}q(\theta)d\theta\right.  \\
+\left.\frac{\varepsilon}{\Gamma(n-\beta)}\left(-\frac{d}{dt}\right)^{n}
\int_t^b (\theta - t)^{n-\beta-1}\xi(\theta,q)d\theta\right]_{\varepsilon=0} \, .
\end{multline}
Expression \eqref{eq:SP} is equivalent to \eqref{eq:cnsidf}.
\end{proof}

The following definition is useful in order to introduce an
appropriate concept of \emph{fractional conserved quantity}.

\begin{defn}
Given two functions $f$ and $g$ of class $C^{1}$ in the interval
$[a,b]$, we introduce the following notation:
\begin{equation*}
\mathcal{D}_{t}^{\gamma}\left(f,g\right)
= -g \, {_tD_b^\gamma} f + f \, {_aD_t^\gamma} g \, ,
\end{equation*}
where $t \in [a,b]$ and $\gamma \in \mathbb{R}_0^+$.
\end{defn}

\begin{rem}
If $\gamma = 1$, the operator $\mathcal{D}_{t}^{\gamma}$ is reduced to
\begin{equation*}
\mathcal{D}_{t}^{1}\left(f,g\right)
=- g \, {_tD_b^1 f} + f \, {_aD_t^1} g = \dot{f} g + f \dot{g}
= \frac{d}{dt}(f g) \, .
\end{equation*}
In particular, $\mathcal{D}_{t}^{1}\left(f,g\right) =
\mathcal{D}_{t}^{1}\left(g,f\right)$.
\end{rem}

\begin{rem}
In the fractional case ($\gamma \ne 1$), functions $f$ and $g$
do not commute: in general $\mathcal{D}_{t}^{\gamma}\left(f,g\right)
\ne \mathcal{D}_{t}^{\gamma}\left(g,f\right)$.
\end{rem}

\begin{rem}
The linearity of the operators $_aD_t^\gamma$ and $_tD_b^\gamma$
imply the linearity of the operator $\mathcal{D}_{t}^{\gamma}$.
\end{rem}

\begin{defn} [fractional-conserved quantity -- \textrm{cf.} Definition~\ref{def:leicon}]
\label{eq:fcl}
We say that $C_{f}\left(t,q,{_aD_t^\alpha q},{_tD_b^\beta q}\right)$
is \emph{fractional-conserved} if (and only if) it is possible to write
$C_{f}$ in the form
\begin{equation}
\label{eq:somaPrd}
C_{f}\left(t,q,d_l,d_r\right)
= \sum_{i=1}^{m} C_{i}^1\left(t,q,d_l,d_r\right) \cdot C_{i}^2\left(t,q,d_l,d_r\right)
\end{equation}
for some $m \in \mathbb{N}$ and some functions $C_{i}^1$ and $C_{i}^2$, $i = 1,\ldots,m$, where
each pair $C_{i}^1$ and $C_{i}^2$, $i = 1,\ldots,m$, satisfy
\begin{equation}
\label{eq:frac-cons}
\mathcal{D}_{t}^{\gamma_i}\left(C_{i}^{j_i^1}\left(t,q,{_aD_t^\alpha q},
{_tD_b^\beta q}\right),C_{i}^{j_i^2}\left(t,q,{_aD_t^\alpha q},
{_tD_b^\beta q}\right)\right) = 0
\end{equation}
with $\gamma_i \in \{\alpha,\beta\}$, $j_i^1 = 1$ and $j_i^2 = 2$ or
$j_i^1 = 2$ and $j_i^2 = 1$,
along all the fractional Euler-Lagrange extremals
(\textrm{i.e.} along all the solutions of the fractional
Euler-Lagrange equations \eqref{eq:eldf}).
\end{defn}

\begin{rem}
Noether conserved quantities \eqref{eq:TeoNet}
are a sum of products, like the structure \eqref{eq:somaPrd}
we are assuming in Definition~\ref{eq:fcl}. For $\alpha = \beta = 1$
\eqref{eq:frac-cons} is equivalent to the standard definition
$\frac{d}{dt} \left[C_{f}\left(t,q(t),\dot{q}(t),-\dot{q}(t)\right)\right] = 0$.
\end{rem}

\begin{exmp}
Let $C_{f}$ be a fractional-conserved quantity written in the form
\eqref{eq:somaPrd} with $m=1$, that is, let $C_{f} = C_{1}^1 \cdot
C_{1}^2$ be a fractional-conserved quantity for some given
functions $C_{1}^1$ and $C_{1}^2$. Condition \eqref{eq:frac-cons}
of Definition~\ref{eq:fcl} means one of four things: or
$\mathcal{D}_{t}^{\alpha}\left(C_{1}^{1},C_{1}^{2}\right) = 0$, or
$\mathcal{D}_{t}^{\alpha}\left(C_{1}^{2},C_{1}^{1}\right) = 0$, or
$\mathcal{D}_{t}^{\beta}\left(C_{1}^{1},C_{1}^{2}\right) = 0$, or
$\mathcal{D}_{t}^{\beta}\left(C_{1}^{2},C_{1}^{1}\right) = 0$.
\end{exmp}

\begin{rem}
Given a fractional-conserved quantity $C_{f}$, the definition of
$C_{i}^1$ and $C_{i}^2$, $i = 1,\ldots,m$, is never unique. In
particular, one can always choose $C_{i}^1$ to be $C_{i}^2$ and
$C_{i}^2$ to be $C_{i}^1$. Definition~\ref{eq:fcl} is imune to the
arbitrariness in defining the $C_{i}^j$, $i = 1,\ldots,m$, $j = 1,2$.
\end{rem}

\begin{rem}
Due to the simple fact that the same function can be written in several
different but equivalent ways, to a given fractional-conserved
quantity $C_{f}$ it corresponds an integer value $m$ in
\eqref{eq:somaPrd} which is, in general, also not unique (see
Example~\ref{ex:dif-m}).
\end{rem}

\begin{exmp}
\label{ex:dif-m}
Let $f$, $g$ and $h$ be functions satisfying
$\mathcal{D}_{t}^{\alpha}\left(g,f\right) = 0$,
$\mathcal{D}_{t}^{\alpha}\left(f,g\right) \ne 0$,
$\mathcal{D}_{t}^{\alpha}\left(h,f\right) = 0$,
$\mathcal{D}_{t}^{\alpha}\left(f,h\right) \ne 0$
along all the fractional Euler-Lagrange extremals
of a given fractional variational problem. One can provide
different proofs to the fact that $C = f\left(g+h\right)$ is a
fractional-conserved quantity: (i) $C$ is fractional-conserved
because we can write $C$ in the form \eqref{eq:somaPrd} with
$m=2$, $C_1^1 = g$, $C_1^2 = f$, $C_2^1 = h$, and $C_2^2 = f$,
satisfying \eqref{eq:frac-cons}
($\mathcal{D}_{t}^{\alpha}\left(C_1^1,C_1^2\right) = 0$ and
$\mathcal{D}_{t}^{\alpha}\left(C_2^1,C_2^2\right) = 0$); (ii) $C$
is fractional-conserved because we can wrire $C$ in the form
\eqref{eq:somaPrd} with $m=1$, $C_1^1 = g+h$, and $C_1^2 = f$,
satisfying \eqref{eq:frac-cons}
($\mathcal{D}_{t}^{\alpha}\left(C_1^1,C_1^2\right)
= \mathcal{D}_{t}^{\alpha}\left(g+h,f\right)
= \mathcal{D}_{t}^{\alpha}\left(g,f\right)
+ \mathcal{D}_{t}^{\alpha}\left(h,f\right) = 0$).
\end{exmp}

\begin{thm}[\textrm{cf.} Theorem~\ref{theo:tnaut}]
\label{theo:tnadf1}
If the functional \eqref{Pf} is invariant
under the transformations \eqref{eq:tinf},
in the sense of Definition~\ref{def:inv1:MR}, then
\begin{multline}
\label{eq:tnaf}
C_{f}\left(t,q,{_aD_t^\alpha q}, {_tD_b^\beta q}\right) \\
= \left[\partial_{3} L\left(t,q,{_aD_t^\alpha q},{_tD_b^\beta q}\right)
- \partial_{4} L\left(t,q,{_aD_t^\alpha q},{_tD_b^\beta q}\right)\right]
\cdot \xi(t,q)
\end{multline}
is fractional-conserved.
\end{thm}

\begin{rem}
In the particular case $\alpha = \beta = 1$, we obtain from
\eqref{eq:tnaf} the conserved quantity \eqref{eq:tnau}
applied to $\mathcal{L}$ \eqref{eq:mathcalL}.
\end{rem}

\begin{proof}
We use the fractional Euler-Lagrange equations
\begin{multline*}
\partial_{2} L\left(t,q,{_aD_t^\alpha q},{_tD_b^\beta q}\right)
= -{_{t}D_b^\alpha}\partial_{3} L\left(t,q,{_aD_t^\alpha q},{_tD_b^\beta q}\right) \\
- {_aD_t^\beta}\partial_{4} L\left(t,q,{_aD_t^\alpha q},{_tD_b^\beta q}\right)
\end{multline*}
in \eqref{eq:cnsidf}, obtaining
\begin{multline*}
-{_{t}D_b^\alpha}\partial_{3} L\left(t,q,{_aD_t^\alpha
q},{_tD_b^\beta q}\right)\cdot\xi(t,q)
+\partial_{3} L\left(t,q,{_aD_t^\alpha
q},{_tD_b^\beta q}\right)\cdot{_aD_t^\alpha \xi(t,q)}\\
- {_aD_t^\beta}\partial_{4} L\left(t,q,{_aD_t^\alpha
q},{_tD_b^\beta q}\right)\cdot\xi(t,q)
+ \partial_{4} L\left(t,q,{_aD_t^\alpha q},{_tD_b^\beta q}\right)
\cdot {_tD_b^\beta}\xi(t,q)\\
= \mathcal{D}_{t}^{\alpha}\left(\partial_{3} L\left(t,q,{_aD_t^\alpha q},
{_tD_b^\beta q}\right),\xi(t,q)\right) 
- \mathcal{D}_{t}^{\beta}\left(\xi(t,q),\partial_{4} L\left(t,q,{_aD_t^\alpha q},
{_tD_b^\beta q}\right)\right) = 0 \, .
\end{multline*}
The proof is complete.
\end{proof}

\begin{defn}[invariance of (\ref{Pf}) -- \textrm{cf.} Definition~\ref{def:inva}]
\label{def:invadf}
Functional \eqref{Pf} is said to be invariant
under the one-parameter group of infinitesimal transformations
\eqref{eq:tinf2} if, and only if,
\begin{multline*}
\int_{t_{a}}^{t_{b}} L\left(t,q(t),{_aD_t^\alpha q(t)},
{_tD_b^\beta q(t)}\right) dt
= \int_{\bar{t}(t_a)}^{\bar{t}(t_b)} L\left(\bar{t},\bar{q}(\bar{t}),
{_aD_{\bar{t}}^\alpha \bar{q}(\bar{t})},{_{\bar{t}}D_b^\beta \bar{q}}(\bar{t})\right) d\bar{t}
\end{multline*}
for any subinterval $[{t_{a}},{t_{b}}] \subseteq [a,b]$.
\end{defn}

The next theorem provides the extension of Noether's theorem
for Fractional Problems of the Calculus of Variations.

\begin{thm}[fractional Noether's theorem]
\label{theo:tndf}
If functional \eqref{Pf} is invariant,
in the sense of Definition~\ref{def:invadf}, then
\begin{multline}
\label{eq:tndf}
C_{f}\left(t,q,{_aD_t^\alpha q},{_tD_b^\beta q}\right)
= \left[\partial_{3} L\left(t,q,{_aD_t^\alpha q},{_tD_b^\beta q}\right) \right. \\
- \left. \partial_{4} L\left(t,q,{_aD_t^\alpha q},
{_tD_b^\beta q}\right)\right] \cdot \xi(t,q) \\
+ \left[L\left(t,q,{_aD_t^\alpha q},{_tD_b^\beta q}\right)
- \alpha\partial_{3} L\left(t,q,{_aD_t^\alpha q},
{_tD_b^\beta q}\right)\cdot{_aD_t^\alpha q}\right. \\
\left.-\beta\partial_{4} L\left(t,q,{_aD_t^\alpha q},{_tD_b^\beta q}\right)
\cdot {_{t}D_{b}^{\beta}q} \right] \tau(t,q)
\end{multline}
is fractional-conserved (Definition~\ref{eq:fcl}).
\end{thm}

\begin{rem}
If $\alpha = \beta = 1$, the fractional conserved quantity
\eqref{eq:tndf} gives \eqref{eq:TeoNet} applied to
$\mathcal{L}$ \eqref{eq:mathcalL}.
\end{rem}

\begin{proof}
Our proof is an extension
of the method used in the proof of Theorem~\ref{theo:tnoe}.
For that we reparameterize the time (the independent variable $t$)
by the Lipschitzian transformation
\begin{equation*}
[a,b]\ni t\longmapsto \sigma f(\lambda) \in [\sigma_{a},\sigma_{b}]
\end{equation*}
that satisfies $t_{\sigma}^{'} = f(\lambda) = 1$ if $\lambda=0$.
Functional \eqref{Pf} is reduced, in this way,
to an autonomous functional:
\begin{multline}
\label{eq:tempo}
\bar{I}[t(\cdot),q(t(\cdot))] =
\int_{\sigma_{a}}^{\sigma_{b}} \hspace*{-0.2cm}
L\left(t(\sigma),q(t(\sigma)),
{_{\sigma_{a}}D_{t(\sigma)}^{\alpha}q(t(\sigma))},
{_{t(\sigma)}D_{\sigma_{b}}^{\beta}q(t(\sigma))}\right)t_{\sigma}^{'} d\sigma ,
\end{multline}
where $t(\sigma_{a}) = a$, $t(\sigma_{b}) = b$,
\begin{equation*}
\begin{split}
_{\sigma_{a}}&D_{t(\sigma)}^{\alpha}q(t(\sigma))\\
&=
\frac{1}{\Gamma(n-\alpha)}\left(\frac{d}{dt(\sigma)}\right)^{n}
\int_{\frac{a}{f(\lambda)}}^{\sigma f(\lambda)}\left({\sigma
f(\lambda)}-\theta\right)^{n-\alpha-1}q\left(\theta f^{-1}(\lambda)\right)d\theta    \\
&=
\frac{(t_{\sigma}^{'})^{-\alpha}}{\Gamma(n-\alpha)}
\left(\frac{d}{d\sigma}\right)^{n}
\int_{\frac{a}{(t_{\sigma}^{'})^{2}}}^{\sigma}
(\sigma-s)^{n-\alpha-1}q(s)ds  \\
&=
(t_{\sigma}^{'})^{-\alpha}{_{\frac{a}{(t_{\sigma}^{'})^{2}}}
D_{\sigma}^{\alpha}q(\sigma)} \, ,
\end{split}
\end{equation*}
and, using the same reasoning,
\begin{equation*}
_{t(\sigma)}D_{\sigma_{b}}^{\beta}q(t(\sigma))=
(t_{\sigma}^{'})^{-\beta}{_{\sigma}
D_{\frac{b}{(t_{\sigma}^{'})^{2}}}^{\beta}q(\sigma)} \, .
\end{equation*}
We have
\begin{equation*}
\begin{split}
&\bar{I}[t(\cdot),q(t(\cdot))] \\
&= \int_{\sigma_{a}}^{\sigma_{b}} \hspace*{-0.3cm}
L\left(t(\sigma),q(t(\sigma)),(t_{\sigma}^{'})^{-\alpha}{_{\frac{a}{(t_{\sigma}^{'})^{2}}}
D_{\sigma}^{\alpha}}q(\sigma),(t_{\sigma}^{'})^{-\beta}{_{\sigma}
D_{\frac{b}{(t_{\sigma}^{'})^{2}}}^{\beta}}q(\sigma)\right) t_{\sigma}^{'} d\sigma \\
&\doteq \int_{\sigma_{a}}^{\sigma_{b}}
\bar{L}_{f}\left(t(\sigma),q(t(\sigma)),t_{\sigma}^{'},{_{\frac{a}{(t_{\sigma}^{'})^{2}}}
D_\sigma^\alpha} q(t(\sigma)),{_\sigma D_{\frac{b}{(t_{\sigma}^{'})^{2}}}^\beta}
q(t(\sigma))\right)d\sigma \\
&= \int_a^b L\left(t,q(t),{_aD_t^\alpha} q(t),{_tD_b^\beta} q(t)\right) dt \\
&= I[q(\cdot)] \, .
\end{split}
\end{equation*}
By hypothesis, functional \eqref{eq:tempo} is invariant
under transformations \eqref{eq:tinf2},
and it follows from Theorem~\ref{theo:tnadf1} that
\begin{multline}
\label{eq:tnadf2}
C_{f}\left(t(\sigma),q(t(\sigma)),t_{\sigma}^{'},{_{\frac{a}{(t_{\sigma}^{'})^{2}}}
D_\sigma^\alpha} q(t(\sigma)),{_\sigma D_{\frac{b}{(t_{\sigma}^{'})^{2}}}^\beta}
q(t(\sigma))\right)\\
= \left(\partial_{4}\bar{L}_{f}-\partial_{5}\bar{L}_{f}\right)
\cdot \xi + \frac{\partial}{\partial t'_\sigma} \bar{L}_{f}\tau
\end{multline}
is a fractional conserved quantity. For $\lambda = 0$,
\begin{equation*}
\begin{split}
_{\frac{a}{(t_{\sigma}^{'})^{2}}}D_\sigma^\alpha q(t(\sigma))
& = {_aD_t}^\alpha q(t) \, ,\\
_\sigma D_{\frac{b}{(t_{\sigma}^{'})^{2}}}^\beta q(t(\sigma))
&= {_tD_b}^\beta q(t) \, ,
\end{split}
\end{equation*}
and we get
\begin{equation}
\label{eq:prfMR:q1}
\partial_{4}\bar{L}_{f}-\partial_{5}\bar{L}_{f}
=\partial_{3} L - \partial_{4} L \, ,
\end{equation}
and
\begin{equation}
\label{eq:prfMR:q2}
\begin{split}
&\frac{\partial}{\partial t'_\sigma} \bar{L}_{f}
= \partial_{4}{\bar{L}_{f}}
\cdot \frac{\partial}{\partial t_{\sigma}^{'}}\left[
\frac{(t_{\sigma}^{'})^{-\alpha}}{\Gamma(n-\alpha)}\left(\frac{d}{d\sigma}\right)^{n}
\int_{\frac{a}{(t_{\sigma}^{'})^{2}}}^{\sigma}
(\sigma-s)^{n-\alpha-1}q(s)ds\right]t_{\sigma}^{'} \\
&+ \partial_{5}{\bar{L}_{f}}\cdot\frac{\partial}{\partial t_{\sigma}^{'}}
\left[\frac{(t_{\sigma}^{'})^{-\beta}}{\Gamma(n-\beta)}\left(-\frac{d}{d\sigma}\right)^{n}
\int_\sigma^{\frac{b}{(t_{\sigma}^{'})^{2}}}
(s - \sigma)^{n-\beta-1}q(s)ds\right]t_{\sigma}^{'} + L \\
&= \partial_{4}\bar{L}_{f}\cdot\left[\frac{-\alpha(t_{\sigma}^{'})^{-\alpha-1}}
{\Gamma(n-\alpha)}\left(\frac{d}{d\sigma}\right)^{n}
\int_{\frac{a}{(t_{\sigma}^{'})^{2}}}^{\sigma}
(\sigma-s)^{n-\alpha-1}q(s)ds\right] \\
&+ \partial_{5}\bar{L}_{f}\cdot
\left[\frac{-\beta(t_{\sigma}^{'})^{-\beta-1}}{\Gamma(n-\beta)}\left(-\frac{d}{d\sigma}\right)^{n}
\int_\sigma^{\frac{b}{(t_{\sigma}^{'})^{2}}}
(s - \sigma)^{n-\beta-1}q(s)ds\right] + L \\
&= -\alpha\partial_{3} L\cdot{_{a}D_{t}^{\alpha}}q
-\beta\partial_{4} L\cdot{_{t}D_{b}^{\beta}}q + L \, .
\end{split}
\end{equation}
Substituting \eqref{eq:prfMR:q1} and \eqref{eq:prfMR:q2} into
equation \eqref{eq:tnadf2}, we obtain
the fractional-conserved quantity \eqref{eq:tndf}.
\end{proof}


\section{Examples}
\label{sec:Ex}

In order to illustrate our results, we consider in this section two
problems from \cite[\S 3]{CD:BalAv:2004}, with
Lagrangians which are linear functions of the velocities.
In both examples, we have used \cite{CD:GouveiaTorres:2005}
to compute the symmetries \eqref{eq:tinf2}.


\begin{exmp}
Let us consider the following fractional problem
of the calculus of variations \eqref{Pf} with $n = 3$:
\begin{multline*}
I[q(\cdot)] = I[q_1(\cdot),q_2(\cdot),q_3(\cdot)] \\
= \int_a^b \left(\left({_{a}D_{t}^{\alpha}q_{1}}\right)q_{2}
-\left({_{a}D_{t}^{\alpha}q_{2}}\right)q_{1}
-(q_{1}-q_{2})q_{3}\right) dt \longrightarrow \min \, .
\end{multline*}
The problem is invariant under the transformations
\eqref{eq:tinf2} with
\begin{equation*}
\left(\tau,\xi_{1},\xi_{2},\xi_{3}\right) = (-ct,0,0,c q_{3}) \, ,
\end{equation*}
where $c$ is an arbitrary constant. We obtain from
our fractional Noether's theorem the following
fractional-conserved quantity \eqref{eq:tndf}:
\begin{multline}
\label{eq:ex1}
C_{f}\left(t,q,{_aD_t^\alpha q}\right)
= \left[(1-\alpha)(\left({_{a}D_{t}^{\alpha}q_{1}}\right)q_{2}
-\left({_{a}D_{t}^{\alpha}q_{2}}\right)q_{1})
-(q_{1}-q_{2})q_{3}\right] t \, .
\end{multline}
We remark that the fractional-conserved quantity \eqref{eq:ex1}
depends only on the fractional derivatives of $q_{1}$ and $q_{2}$
if $\alpha \in ]0,1[$. For $\alpha=1$, we obtain the classical result:
\begin{equation*}
C(t,q)=\left(q_{1}-q_{2}\right) q_{3} t
\end{equation*}
is conserved for the problem of the calculus of variations
\begin{equation*}
\int_a^b \left(
\dot{q}_{1}q_{2}-\dot{q}_{2}q_{1}-(q_{1}-q_{2})q_{3} \right) dt
 \longrightarrow \min \, .
\end{equation*}
\end{exmp}


\begin{exmp}
We consider now a variational functional \eqref{Pf} with $n = 4$:
\begin{multline*}
I[q(\cdot)] = I[q_1(\cdot),q_2(\cdot),q_3(\cdot),q_4(\cdot)] \\
=- \int_a^b \left[\left({_{t}D_{b}^{\beta}q_{1}}\right)q_{2}
+\left({_{t}D_{b}^{\beta}q_{3}}\right)q_{4}
-\frac{1}{2}\left(q_{4}^{2}-2q_{2}q_{3}\right)\right]dt \, .
\end{multline*}
The problem is invariant under \eqref{eq:tinf2} with
\begin{equation*}
\left(\tau,\xi_{1},\xi_{2},\xi_{3},\xi_{4}\right)
= \left(\frac{2c}{3}t,cq_{1},-cq_{2},\frac{c}{3}q_{3},-\frac{c}{3}q_{4}\right) \, ,
\end{equation*}
where $c$ is an arbitrary constant.
We conclude from \eqref{eq:tndf} that
\begin{multline}
\label{eq:ex2}
C_{f}\left(t,q,{_tD_b^\beta q}\right)
=\Big[(\beta-1)\left(\left({_{t}D_{b}^{\beta}q_{1}}\right)q_{2}
+ \left({_{t}D_{b}^{\beta}q_{3}}\right)q_{4}\right) \\
+ \frac{1}{2}\left(q_{4}^{2}-2q_{2}q_{3}\right)\Big]\frac{2}{3}t
+ q_{1}q_{2}+\frac{q_{3}q_{4}}{3}
\end{multline}
is a fractional conserved quantity. In the particular case $\beta=1$,
the classical result follows from \eqref{eq:ex2}:
\begin{equation*}
C(t,q)= q_{1}q_{2}+\frac{q_{3}q_{4}}{3}
+ \frac{1}{3}\left(q_{4}^{2}-2q_{2}q_{3}\right) t
\end{equation*}
is preserved along all the solutions of the Euler-Lagrange
differential equations \eqref{eq:el} of the problem
\begin{equation*}
\int_a^b \dot{q}_{1}q_{2}+\dot{q}_{3}q_{4}
+ \frac{1}{2}\left(q_{4}^{2}-2q_{2}q_{3}\right) dt
\longrightarrow \min \, .
\end{equation*}
\end{exmp}


\section{Conclusions and Open Questions}
\label{sec:Conc}

The fractional calculus is a mathematical area of a currently strong
research, with numerous applications in physics and engineering.
The theory of the calculus of variations for fractional systems
was recently initiated in \cite{CD:Agrawal:2002}, with the proof
of the fractional Euler-Lagrange equations. In this paper we go a
step further: we prove a fractional Noether's theorem.

The fractional variational theory is in its childhood so that
much remains to be done. This is particularly true
in the area of fractional optimal control,
where the results are rare. A fractional
Hamiltonian formulation is obtained in \cite{Muslih:2005},
but only for systems with linear velocities.
A study of fractional optimal control problems
with quadratic functionals can be found in \cite{Agrawal:2004a}.
To the best of the author's knowledge, there is no
general formulation of a fractional version
of Pontryagin's Maximum Principle. Then, with a
fractional notion of Pontryagin extremal,
one can try to use the techniques of \cite{delfimEJC}
to extend the present results to the more
general context of the fractional optimal control.


\subsection*{Acknowledgements}

GF acknowledges the financial support of IPAD
(Portuguese Institute for Development);
DT the support from the Control Theory Group
(cotg) of the Centre for Research
in Optimization and Control (http://ceoc.mat.ua.pt).
The authors are grateful to Paulo Gouveia and Ilona Dzenite
for several suggestions.



\end{document}